\pdfoutput=1

\documentclass[a4,notitlepage,12pt]{jedm}
\usepackage[utf8]{inputenc}
\usepackage{pdfpages}
\usepackage{siunitx}
\usepackage{graphicx}
\usepackage{amsmath}
\usepackage{amssymb}
\usepackage{booktabs}
\usepackage[export]{adjustbox}
\usepackage{subcaption}
\usepackage{tabu}
\usepackage[T1]{fontenc}
\usepackage{placeins}
\usepackage[labelfont={bf}, margin=1cm]{caption}
\usepackage{url}
\usepackage{adjustbox}
\usepackage[colorlinks = true, linkcolor = blue, urlcolor  = blue, citecolor = magenta]{hyperref}
\usepackage{textcomp}
\usepackage{csvsimple, array}
\usepackage{floatrow}
\usepackage{blindtext}
\usepackage[noadjust]{cite}
\title{Core Course Analysis for Undergraduate Students in Mathematics}

\author{{\large Ritvik Kharkar}\\University of California, Los Angeles\\ritvikmath@ucla.edu \and {\large Jessica Tran}\\University of California, Los Angeles\\jessicatran@ucla.edu  \and {\large Charles Z. Marshak}\\University of California, Los Angeles\\cmarshak@math.ucla.edu}

\date{\today}

\renewcommand{\t}[1]{\mathrm{#1}}
\renewcommand{\vec}[1]{\mathbf{#1}}
\DeclareMathOperator*{\argmin}{arg\,min}

\begin{document}

\maketitle

\begin{abstract}
In this work, we develop statistical tools to understand core courses at the university level. Traditionally, professors and administrators label courses as "core" when the courses contain foundational material.  Such courses are often required to complete a major, and, in some cases, allocated additional educational resources.  
We identify two key attributes which we expect core courses to have. Namely, we expect core courses to be highly correlated with and highly impactful on a student's overall mathematics GPA. We use two statistical procedures to measure the strength of these attributes across courses.
The first of these procedures fashions a metric out of standard correlation measures. The second utilizes sparse regression. 
We apply these methods on student data coming from the University of California, Los Angeles (UCLA) department of mathematics to compare core and non-core coursework. 

\end{abstract}

\section{Introduction}
Educators, professors, and mathematicians are continually improving mathematics education  \cite{ellis2005paradigm, wu1997mathematics}. There is a consensus that mathematics education divides into two categories: design of math curriculum and the development of teaching strategies.  In particular, members from the Mathematical Association of America and the American Mathematical Society \cite{roberts2001eh} are regularly analyzing the core requirements of undergraduate mathematics and constantly improving the effectiveness of mathematics teachers \cite{wu1999education, thames2013making, edwards1994current}. In this paper we focus on the design of the curriculum using tools to understand core courses using grade statistics.

The undergraduate math curriculum is the set of courses required for an undergraduate to earn his or her degree.  Departmental faculty and counselors carefully select courses which provide a well-rounded mathematics education. Traditionally, core courses within the mathematics curriculum are those that contain the most foundational material and that inform subsequent coursework.  These courses are often ones that the department encourages students to take early to ensure they get the most from their major.  

 We will analyze the core curriculum using techniques and tools in Educational Data Mining (EDM). EDM is an interdisciplinary research area that aims to improve education utilizing data analysis \cite{romero2010educational, romero2007educational, xu2005exploration, liu2014fitting}. Although researchers have been collecting and analyzing educational data for a long time, it has only recently been established as its own field, due to the accessibility of statistical tools and the emergence of EDM journals and conferences \cite{ baker2009state, siemens2012learning}.
 Tools developed in EDM research have been more quickly adopted due to the challenge of serving a growing undergraduate population \cite{nsf_college}. Recent advancements in education made possible through data analysis tools include course recommendation systems \cite{farzan2006social, xu2015personalized, chu2003designing}, student clustering \cite{mcgill1978variations}, grade predictions \cite{bridgeman2000predictions, meier2015predicting}, and personalized education through online medium \cite{ tekin2015etutor, weld2012personalized}. Adapting such data mining tools for assessing student performance is a difficult task but will help administrators even better serve their student body. 

The focus of this paper will be to provide statistical criteria that validates the department's core course selection. Specifically, we compared the grade distribution from a single math course to the distribution of GPAs. We applied this methodology to the grades of UCLA mathematics students from 2000-2015.  The viewpoint we develop here can be easily adapted to any other department for a similar analysis. 
%

To compare course statistics and GPA statistics we employ two fundamental statistical tools: the correlation coefficient and linear regression. The correlation coefficient is a way to measure the covariance of two distributions \cite{hinkle}.  Regression is a  tool that linearly models how particular measured inputs influence a response \cite{hinkle}.  It informs some of our most important models for business, biology, machine learning, and environmental science \cite{dimitras1996survey, henderson2007application, tu1996advantages}.  Ultimately, we will use these tools to compare the grades earned by students in a single class with these students' overall GPAs.

The data-centric viewpoint we develop offers a statistical methodology to validate the core courses within a university's curriculum. Using this methodology, we can identify those courses that especially impact a student's study. Any discrepancy between the core courses defined by the department and those found from this statistical viewpoint can prompt further discussion of an effective way to design the curriculum. Ultimately, by analyzing the university curriculum, we hope to aid departments in better meeting their goals and preparing exiting students better for industry and academic jobs.

We now discuss how this paper is organized.  In section \ref{detect} of our paper, we describe our data-centric approach to core course analysis providing justification for the use of correlation coefficients and regression. In section \ref{detect}, we precisely describe the statistical setup and provide some sample output.  In section \ref{core_courses}, we apply this methodology of determining core courses to pure mathematics majors. 







\section{Characterizing Core Courses}
\label{detect}
Core courses are those that are the most foundational for a well-rounded mathematics education and are often allocated additional educational resources. 
In this study, we consider the core courses selected by the UCLA  mathematics department and identify statistical attributes that we expect such courses to have. 
We develop statistical tools to measure the strength of these attributes across courses. Using these tools, we then compare the core and non-core courses to validate the department's core course selection. At UCLA, the core courses are Linear Algebra and Real Analysis. All mathematics students must take these courses and are recommended to take them early in their course sequence. 



The first attribute we consider is the correlation between the course grade of a particular student and his or her overall mathematics GPA. By investigating this correlation, we analyze the relationship between rankings determined by grade compared to rankings determined by overall mathematics GPA. 
We expect that the course grades of core courses defined by the department would be highly correlated with one's overall mathematics GPA.

The second attribute we consider is the impact of a particular course grade on overall mathematics GPA. The impact of a course grade on a student's overall mathematics GPA signifies how sensitive overall mathematics GPA is to the course grade.  We examine impact by performing lasso regression in order to select the most important courses.  When we input courses from the mathematics curriculum as variables into lasso, the algorithm selects a subset of courses which are most relevant to overall mathematics GPA. We conclude that the courses assigned the highest coefficients by lasso regression are those having the greatest bearing on overall mathematics GPA. We now describe the correlation metric used to measure the first attribute.


\subsection{Mixed Correlation Metric}
\label{mixed_metric}
To measure the correlation between course grades and a student's overall mathematics GPA, we combine two correlation coefficients. We then analyze these correlation coefficients to measure how well the ranking determined by a single course grade compares to the ranking of overall mathematics GPA. We first consider Spearman's $\rho$ \cite{spearman} which we will denote by $\rho_S$. Given a sample of $n$ students, we can rank students both by overall mathematics GPA as well as by grade in a particular course. Let $r_{1}^{\t {GPA}}, ... , r_{n}^{\t{GPA}}$ be the ranking of students by overall mathematics GPA and let $r^{C}_{1}, ... , r_{n}^{C}$ be the ranking of students by grade in course $C$. The measure $\rho_S$ is defined as
$$
\rho_S= 1-\frac{1}{Z}\left.\sum_{i=1}^{n}(r_{i}^{\t{GPA}} - r_{i}^{C})^2\right.
$$
where $Z = \frac{n(n^2-1)}{6}$ normalizes the sum.  The coefficient $\rho_S$ ranges from $-1$ to $1$ in which positive values indicate an agreement between the two rankings. 
In the context of our study, we use  $\rho_S$ to measure how much a rank determined by a student's grade in a course agrees with his or her overall mathematics GPA. 

We will also employ the coefficient of determination, $R_{sq}$, a normalized measure of covariance between grade in a particular course and overall mathematics GPA. Let $g_{1}^{C}, ... , g_{n}^C$ be the grade in course $C$ of each student and let $g_{1}^{\t{GPA}}, ... , g_{n}^{\t{GPA}}$ be the overall mathematics GPA of each student. The sample means of these quantities are defined as
\begin{align*}\overline{{g}}\;{}^C &= \frac{1}{n}\sum_{i=1}^{n}g_{i}^C\\ 
\overline{g}\;{}^{\t{GPA}} &= \frac{1}{n}\sum_{i=1}^{n}{g_{i}^{\t{GPA}}}.
\end{align*}
$R_{sq}$ is defined as follows, 
$$ 
R_{sq} = \frac{\sum_{i=1}^{n}(g^C_i - \overline{g}\;{}^{C})(g_{i}^{\t{GPA}} - \bar{g}\;{}^{\t{GPA}})}{\sqrt{\sum_{i=1}^{n}\left(g^C_{i}-\bar{g}\;{}^{C}_i\right)^2}\;\sqrt{\sum_{i=1}^{n}\left(g_{i}^{\t{GPA}}-\bar{g}\;{}^{\t{GPA}}\right)^2}}.
$$
Using these two correlation coefficients, we define a new holistic metric as $\rho_H = \sqrt{{(\rho_S)}^2 + (R_{sq})^2}$. This new metric, $\rho_H$ takes into account both correlation metrics and is the $\ell_2$ distance to the origin. 

\begin{figure}
    \centering
    \includegraphics[scale = 0.7]{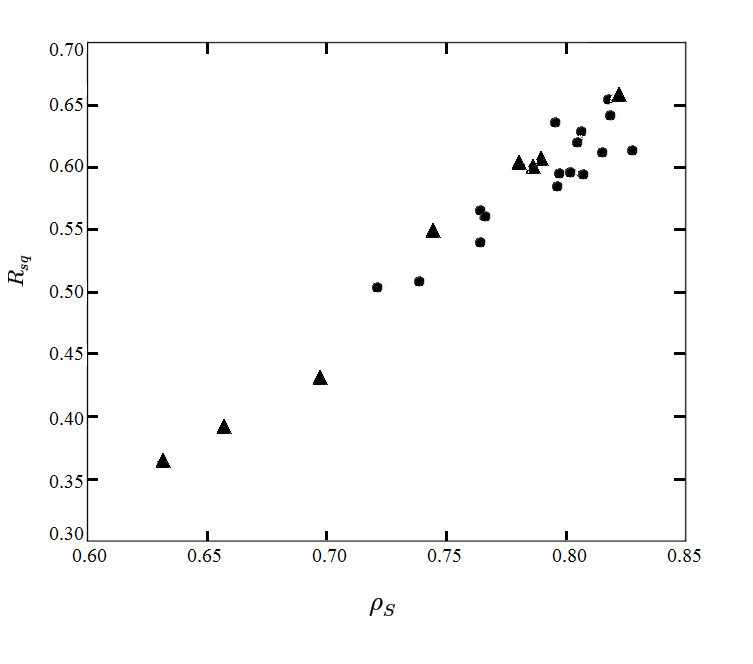}
    \caption{Scatter plot of $\rho_S$ values against $R_{sq}$ values for the Math/Econ, Pure Math, Applied Math courses. The triangular points are courses which are required for all majors and the round points are courses which are not required for all majors.}
    \label{fig:scatter_example}
\end{figure}

We give an example of how these metrics are used to identify core courses. We will only consider students who have taken at least five mathematics courses and at least one upper division mathematics course at UCLA so we can identify those courses that might impact students who have made a significant time commitment to their mathematical development.  We provide an example that includes the three largest mathematics majors at UCLA: Mathematics/Economics (Math/Econ), Pure Mathematics (Pure Math), and Applied Mathematics (Applied Math) and will only consider students admitted as freshmen. We omit transfer students since these students usually have different training than students that are admitted right from high school. We report the findings in Tables \ref{pm_mc1} and \ref{pm_mc2}. After generating the $\rho_S$ and $R_{sq}$ values for each student group in each of the courses, we plot the two metrics to understand how well the rankings matched. The scatter plot in Figure \ref{fig:scatter_example} plots $\rho_S$ against $R_{sq}$ using the values from Tables \ref{pm_mc1} and \ref{pm_mc2}.

\FloatBarrier

\newcommand{\lengthmethodtable}{4.2cm}
\newcommand{\shortlengthmethodtable}{1.4cm}
\begin{table*}
\centering
\caption{Mixed Correlation Coefficients}
\begin{tabular}{m{4.5cm} cccc}
\toprule
     \textbf{Course} &\textbf{$\rho_{S}$} & \textbf{$R_{sq}$} & \textbf{$\rho_{H}$}\\
    \hline
\begin{filecontents*}{mix_corre_ex1.csv}
Course,p_S,R_sq,p_H
Linear Algebra ,0.821,0.655,1.05
Complex Analysis,0.819,0.657,1.049
Ord. Diff. Eq.,0.818,0.654,1.047
Num. Analysis (Part II),0.827,0.614,1.03
Probability (Part I),0.806,0.62,1.02
Probability (Part II),0.814,0.607,1.015
Game Theory,0.805,0.619,1.015
Num. Analysis (Part I),0.794,0.62,1.007
Nonlinear/Linear Diff. Eq.,0.805,0.603,1.006
Optimization,0.801,0.601,1.001
Combinatorics,0.797,0.601,0.998
Diff. Eq.,0.787,0.61,0.996
\end{filecontents*}
\csvreader[ late after line=\\\hline,table foot = \bottomrule]
{mix_corre_ex1.csv}{Course = \Course, p_S = \pS,  R_sq= \Rsq, p_H = \pH }
{ \Course & \pS & \Rsq &\pH }
\end{tabular}
\label{pm_mc1}
\end{table*}

\begin{table*}
\centering
\caption{Mixed Correlation Coefficients}
\begin{tabular}{m{4.5cm} cccc}
\toprule
     \textbf{Course} &\textbf{$\rho_{S}$} & \textbf{$R_{sq}$} & \textbf{$\rho_{H}$}\\
    \hline
\begin{filecontents*}{mix_corre_ex2.csv}
Course,p_S,R_sq,p_H
Diff. Geometry (Part I),0.802,0.582,0.991
Math Modeling,0.797,0.584,0.988
Real Analysis (Part I),0.781,0.6,0.984
Multivar. Calc. (Part II),0.785,0.589,0.981
Algebra for App.,0.765,0.552,0.943
Real Analysis (Part II),0.764,0.546,0.939
Linear Algebra,0.744,0.548,0.924
Discrete Math,0.738,0.506,0.895
Abs. Algebra (Part I),0.721,0.512,0.885
Multivar. Calc. (Part I),0.696,0.423,0.814
Single Var. Calc. (Part II),0.657,0.383,0.76
Single Var. Calc. (Part I),0.629,0.34,0.715
\end{filecontents*}
\csvreader[ late after line=\\\hline,table foot = \bottomrule]
{mix_corre_ex2.csv}{Course = \Course, p_S = \pS,  R_sq= \Rsq, p_H = \pH }
{ \Course & \pS & \Rsq &\pH }
\end{tabular}
\label{pm_mc2}

\end{table*}

\FloatBarrier

\begin{figure}
    \centering
    \includegraphics[scale = 0.6]{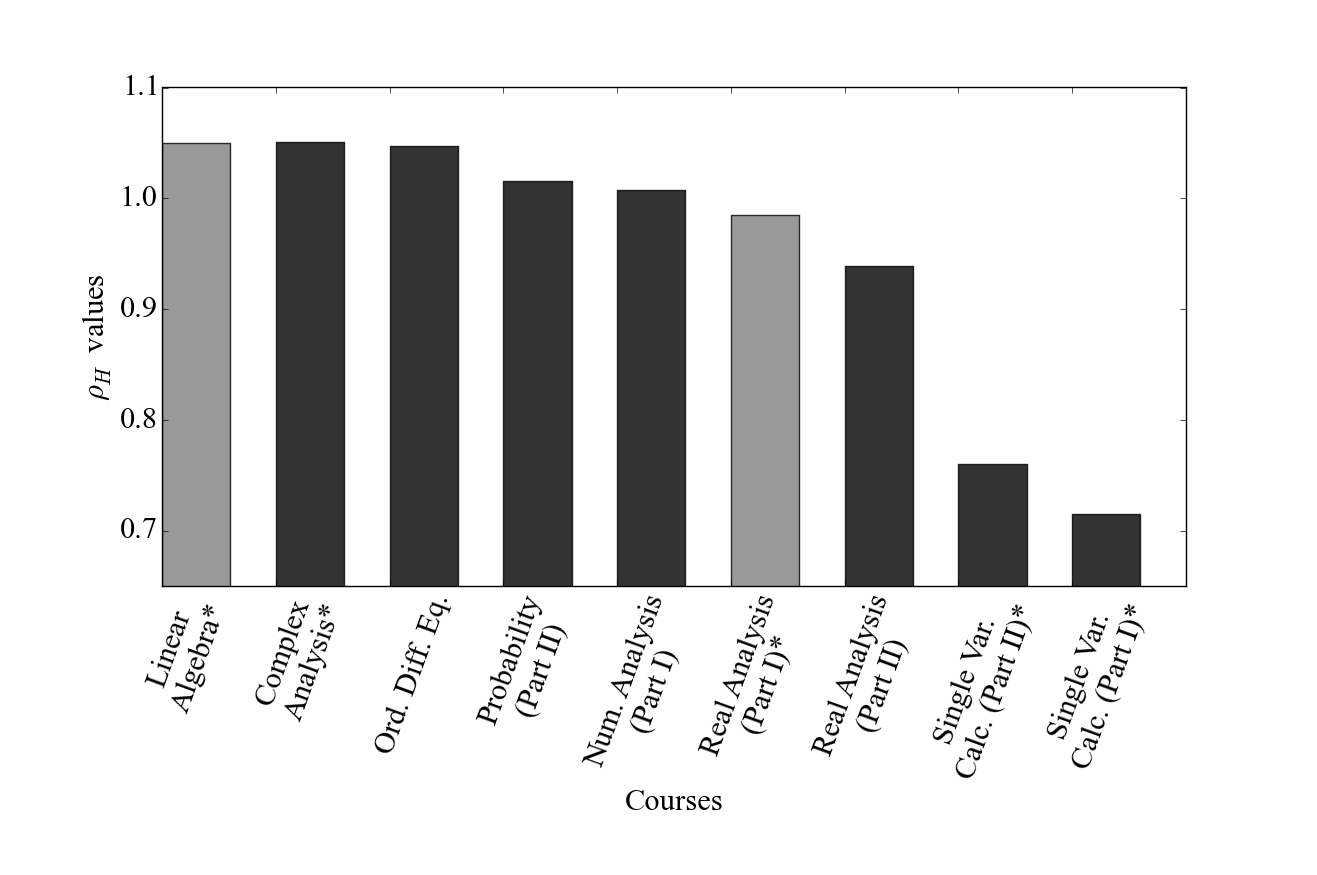}
    \caption{Bar plot of $\rho_{H}$ values of a sample of courses from the correlation metric. The top courses were calculated based on the distance from the origin. Bars colored grey are core courses as defined by the mathematics department. The asterisk (*) indicates courses that are required. 
    }
    \label{fig:bar_plot}
\end{figure}

We can calculate $\rho_{S}$ and $R_{sq}$ for each course. We can then calculate the $\ell_2$ distance from the origin to ($\rho_{S}$,$R_{sq}$) to get a measure of how accurately grades in each course correlate with overall mathematics GPA values. Thus, the larger this distance, the stronger the correlation between course grade and overall mathematics GPA. Figure \ref{fig:bar_plot} shows some of the courses closest to the origin. Since our sample population is from a diverse set of majors, we indicate overlapping required courses with a (*). In the discussion of our results, we use course labellings as assigned by the math department of UCLA, see Appendix \ref{course_appendix}. 


From this bar plot, we see that the only required upper division courses that appear from this metric are  Linear Algebra, Complex Analysis and Real Analysis. This validates that Linear Algebra is a foundational course and also agrees with the department's decision to recommend this course as a first upper division course. Complex Analysis, although not considered a core course by the department, follows immediately after Linear Algebra. On the other hand, Real Analysis, determined a core course by the department, ranks further down in our list. We will perform regression analysis to find which courses are most impactful.

\FloatBarrier
\subsection{Lasso Regression Metric}

\begin{figure}[h!]
\centering
\includegraphics[scale = 0.7]{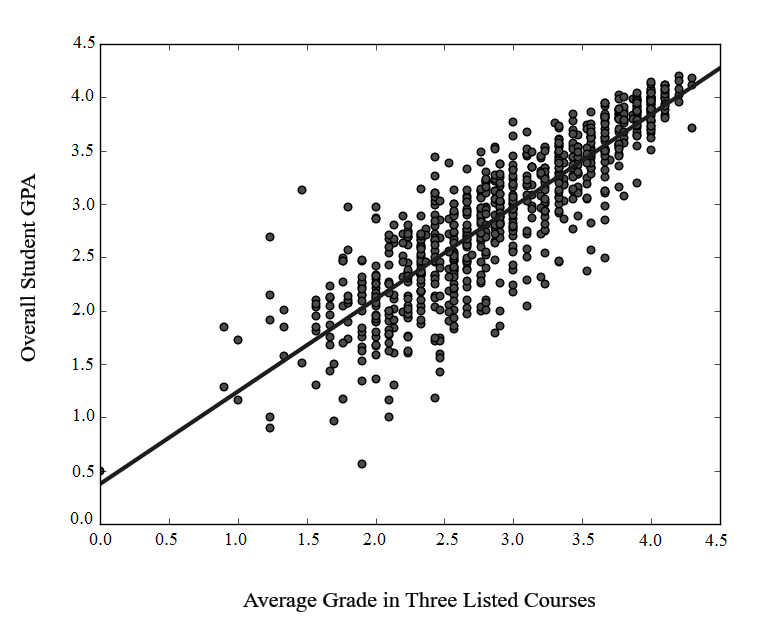}
\caption{Scatter plot of average grade in Linear Algebra, Real Analysis, and Probability Theory (Part I) against overall student GPA for all students who took all three courses. The ridge regression prediction function is drawn through the points.}
\label{fig:ridge}
\end{figure}

In order to further validate the department's selection of core courses, we will now turn our attention to analyzing the impact of a particular course's grade on overall mathematics GPA.
We model the overall mathematics GPA as a linear response to the different course grades a student earns.
We chose to use a linear model due to the versatility of such models, the fact that GPA is a linear function of all course grades, and the apparent linear relationship between core course grades and overall mathematics GPA as seen in Figure \ref{fig:ridge}.
We proceed by selecting the most impactful courses using sparse regression. This regression serves a dual purpose as it also assigns zero coefficient values to courses which have minimal bearing on overall mathematics GPA. This effectively allows us to disqualify certain courses from being considered core. A key application of lasso regression is grade prediction. That is, given the impactful courses selected by lasso regression, we can input the grades in these courses into a different regression model aimed towards predicting overall mathematics GPA. This regression model, better suited to calculate the exact values of the coefficients for each course, will allow us to estimate a student's overall mathematics GPA based on his or her grade in a few courses. 

\FloatBarrier

\subsubsection{Course Selection with Lasso}
\label{Lasso Regression}

We proceed by selecting the most impactful courses with lasso regression \cite{hans2009bayesian}. 
Lasso regression uses an $\ell_1$ penalty to ensure the coefficients for the linear model are sparse.  
A non-zero regression coefficient indicates that the corresponding course is impactful while a zero-coefficient indicates the corresponding course is not.  By removing the courses that have little bearing on a student's overall mathematics GPA, we are able to extract the most impactful courses and validate core course selection by the mathematics department.

Let $C_1,\hdots, C_m$ be the $m$ available courses in the math department.  We can model a student $i$'s overall GPA $g_i^{\t{GPA}}$ as

$$ g^{\t{GPA}}_i = \beta_0 + \beta_{1}g^{C_1}_i + \beta_2g^{C_2}_i + ...  + \beta_mg^{C_m}_i + \epsilon_i $$
where 
$$
g^{C_k}_i = \begin{cases}
$student $i$'s grade in class $k$ $&$ if student $i$ took class $k\\
0 &$ otherwise$
\end{cases}
$$
When considering all $n$ students, we can write this model in the matrix form
$$\vec{g}^{\t{GPA}} = \vec{G}\vec{\beta} + \vec{\epsilon}$$
where the matrix $\vec{G} \in \mathbb{R}^{n\times m}$ contains all of our grade data for all $m$ math courses at UCLA for all $n$ students.


The coefficient vector $\beta$ is determined as
\begin{align}
    \argmin_{\beta} ||\vec{g^{\t{GPA}}} - \vec{G}\vec{\beta}||^2_2 + \alpha||\vec{\beta}||_1
    \label{lasso_min}
\end{align}
for $\alpha>0$. 
For larger $\alpha$, we obtain sparser solutions.  
As $\alpha$ approaches $0$, we obtain solutions that resemble those obtained by ordinary least squares. 
While there is no closed form solution for Eq. \ref{lasso_min}, the minimization can be solved using standard techniques in convex optimization \cite{boyd2011distributed}. 
Throughout this paper, we select $\alpha = 0.1$ when using a lasso model because, experimentally, it has given that approximately 80\% of the coefficients are zero.
We note that the $\beta$ coefficients produced by lasso regression have a smaller variance than those produced by an ordinary least squares regression \cite{tibshirani1996regression}. Unfortunately, this comes at the price of the lasso regression $\beta$'s being more biased estimators.

The results of the lasso regression for freshman mathematics students from the majors Pure Mathematics, Applied Mathematics, and Mathematics/Economics are shown in Table \ref{table:fresh}.


\newcommand{\mediumlengthmethodtable}{3.1cm}
\begin{table*}
\centering
\caption{Lasso regression coefficient values for Pure Mathematics, Mathematics/Economics, and Applied Mathematics  freshman admits. }
\begin{tabular}{m{3.5cm} cc}
\toprule
     \textbf{Course} &\textbf{Lasso Coefficient}\\
    \hline
\begin{filecontents*}{lasso_csv_edit.csv}
Course,Lasso Coefficient
Complex Analysis,0.05734
Diff. Geometry (Part I),0.0502
Num. Analysis (Part II),0.04429
Diff. Eq. (Lower Div.),0.04375
Probability (Part I),0.03322
Math Modeling,0
Multivar. Calc. (Part I),0
Math History,0
\end{filecontents*}
\csvreader[ late after line=\\\hline,table foot = \bottomrule]
{lasso_csv_edit.csv}{Course = \Course, Lasso Coefficient = \LassoCoefficient,  }
{ \Course & \LassoCoefficient }
\end{tabular}

\label{table:fresh}
\end{table*}

\subsubsection{Grade Prediction with Ridge}
\label{Ridge Regression}
Given the courses selected by lasso regression, we will now utilize ridge regression to assign accurate numerical weights to this subset of courses and thereby gauge precisely how impactful each course in our subset is for overall mathematics GPA. In addition, we will also have a meaningful prediction function for overall mathematics GPA based on grades in a few impactful courses.  

We use ridge regression because it will help us account for the multicollinearity between the grades in particular courses, our explanatory variables. We expect multicollinearity to exist due to the fact that students who do well in one course tend to do well in others, and vice versa. This multicollinearity, if unchecked, can skew our regression coefficients, and produce false values \cite{multico}. 

It is important to note that while ridge regression helps to correct for multicollinearity issues by imposing an $\ell_2$ penalty, it causes the resulting regression coefficients to be biased \cite{gaussmarkov}, systematically shifted towards lower values. But, upon running our tests with both ridge regression and ordinary least square regression, we find that ridge regression preserves the relative ordering of the coefficients, produces a higher accuracy level as measured by the coefficient of determination $R_{sq}$, and deviates from the ordinary least squares $\vec{\beta}$ values only slightly.

In our grade prediction process, we will consider three courses as input into our ridge model. The first two courses will always be Linear Algebra (Part I) and Real Analysis (Part I), because they are considered core courses by the mathematics department. We will cycle through each of the top three courses and last course in Table \ref{table:fresh} as our third course in the ridge model, thereby running the ridge model four separate times. We run the model four times because the lasso model selected the top three courses, indicating that overall mathematics GPA is especially sensitive to grades in each selected course. We run the model a fourth time with a course that lasso indicated as non-core to check whether our accuracy level diminishes as expected. We expect our accuracy level to diminish as we pick courses with progressively lower lasso coefficients. 

We use a linear ridge model so that our model for student $i$ is 
$$ g^{\t{GPA}} = \beta_0 + \beta_{1}g^{115A}_i + \beta_2g^{131A}_i + \beta_3g^{C_3}_i + \epsilon_i. $$
where $C_3$ takes on the label of each of the top three courses in Table \ref{table:fresh} as well as the last course in Table \ref{table:fresh}.




We can compactly express the implied system of equations for $n$ students as
$$\vec{g^{\t{GPA}}} = \vec{G}\vec{\beta} + \vec{\epsilon}$$
where $\vec{G} \in \mathbb{R}^{n\times4}$ is a matrix encoding each student's grade in each of the three courses we are using. 





The $\vec{\beta}$ coefficient vector is given by
\begin{align}
    \argmin_{\beta} ||\vec{\vec{g^{\t{GPA}}}} - \vec{G}\vec{\beta}||^2_2 + \alpha||\vec{\beta}||_2
    \label{ridge_min}
\end{align}
for $\alpha>0$. 
We find that $\alpha$ = 2.5 gives the highest accuracy in our predicted overall mathematics GPA values, ${\hat{g}}^{\t{GPA}}$.
For each separate run of the ridge model, we consider all students in the population who have taken Linear Algebra, Real Analysis, and the rotating third course. We then take a random 90\% sample of this group as our training set and use it to generate a coefficient vector $\vec{\beta}$ using Eq. \ref{ridge_min}. The remaining 10\% of the data is used to test this function. We repeat this procedure 100 times, taking a different 90\% random sample of the students each time, and then take the average of the 100 resulting $\vec{\beta}$ vectors to generate an averaged prediction model for overall mathematics GPA. We also average the accuracy levels, measured by $R_{sq}$, over 100 runs of a particular combination of three courses. 

We find that by using Complex Analysis, the course with the highest lasso coefficient in Table  \ref{table:fresh}, as our third course, we get the highest $R_{sq}$ value. This accuracy level drops by 5\% when we use Differential Geometry (Part I), having the second highest lasso coefficient. Furthermore, when we use History of Mathematics as our third course, listed as having a zero lasso coefficient, our accuracy drops by 17\% relative to the model with Complex Analysis.

Our ridge model thus serves to validate the results of our lasso regression. The fact that courses assigned higher lasso coefficients also lead to higher grade prediction accuracy levels provides evidence that our lasso coefficients are accurately ranked. Furthermore, the dramatic decline in accuracy attributed to using courses with zero lasso coefficient in the ridge model verifies that these courses have little bearing on a student's overall mathematics GPA. 
\FloatBarrier

\FloatBarrier

\section{Core Course Analysis for UCLA Math Students}
\label{core_courses}

Now that we have established  methods to validate the core courses, we apply these methods on several different math majors at UCLA. The UCLA mathematics department houses seven different majors. While each major has a different set of course requirements, the department considers Real Analysis and Linear Algebra as the core courses for each major and these two courses are required for \emph{all} majors. We expect grades in these two core courses to be highly correlated with and highly impactful of a student's overall mathematics GPA, more so than non-core courses. We compare the two core and the other non-core courses within a major using the statistical tools outlined in the previous section.  We highlight those departmental non-core courses that have high $\beta_i$ and $\rho_H$ values relative to departmental core courses and interpret this apparent discrepancy.

\subsection{Core Course Analysis for Pure Mathematics Freshmen}\label{pm_study}

We now  analyze core and non-core courses for those students majoring in pure mathematics that were admitted to UCLA as freshmen. We found that our course analysis was the most interesting for these students because many non-core courses were both highly correlated to and highly impactful of these students' overall mathematics GPA, more so than the two departmental core courses. We carefully discuss the findings for this major below.  We also provide some possible explanations as to why many non-core courses exhibited higher correlation and impact on overall mathematics GPA than core courses did.  



\begin{figure}
    \centering
    \includegraphics[scale = 0.6]{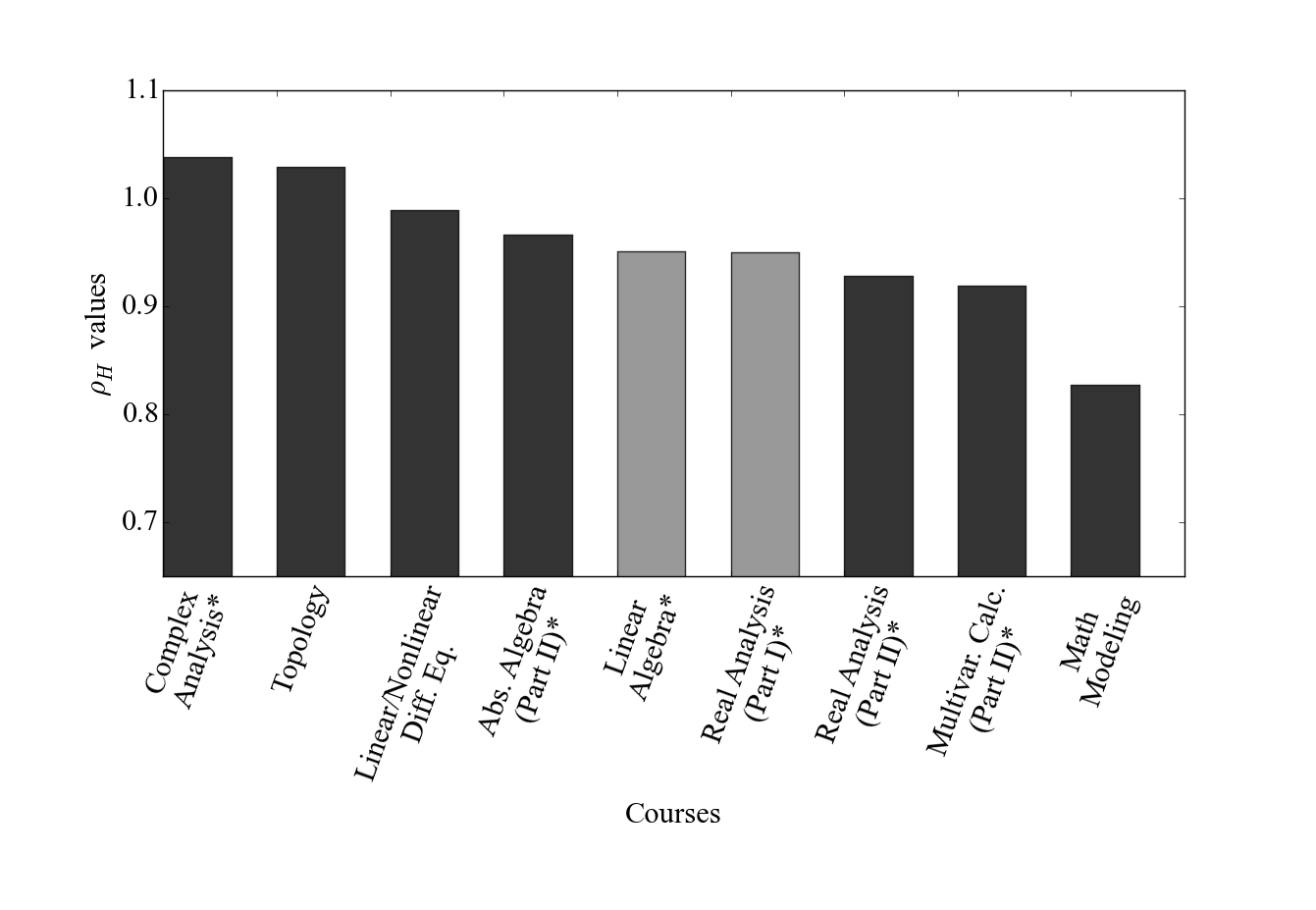}
    \caption{Bar plot of courses of pure mathematics freshman student admits.  Bars colored grey are core courses as defined by the mathematics department. Courses marked with an asterisk (*) are required courses.  }
    \label{fig:pm_ov}
\end{figure}

First, we study $\rho_H$ of courses when considering pure mathematics students.  We display the $\rho_{H}$-values in Figure \ref{fig:pm_ov}. In this bar plot, we show the nine courses with the greatest $\rho_H$.  We see that the two core courses Linear Algebra (Part I) and Real Analysis (Part I) are in this list, but are not the top ranked.  Surprisingly, the top ranked courses are Complex Analysis, Topology, and Linear/Nonlinear Differential Equations.  While every pure math student is required to take Complex Analysis, they are not required to take it early in their major as it is not required for later coursework.  Topology and Linear/Nonlinear Differential Equations are not required at all.  While core courses were highly correlated with a student's overall mathematics GPA, we were surprised that Complex Analysis scored highest. Unlike most other courses in the pure mathematics major, Complex Analysis has no upper-division prerequisite courses.  Indeed, until recently, it was also required for Electrical Engineering students who generally do not take the two core mathematics courses. Moreover, no other upper division mathematics course at UCLA requires Complex Analysis.  On the surface, we expect that Complex Analysis could be taken in a student's last quarter without impacting their performance in other courses.  As such, we also expect that such a course would have lower correlation to overall mathematics GPA.  

Next, we study the lasso coefficients $\beta_i$ and display the results in Table \ref{table:results_final}. Here, the course with the highest $\beta_{i}$ value is again Complex Analysis, followed by Probability Theory (Part I), Honors Abstract Algebra (Part I), and Differential Equations. Linear Algebra and Real Analysis, the core courses, both have zero as their coefficients, suggesting that non-core courses have greater impact on overall mathematics GPA. More formally, those courses with non-zero coefficients are the most important inputs for a linear model of overall mathematics GPA. We are surprised to find that Complex Analysis ranks as the top course under both the $\rho_H$ and $\beta_i$ metrics because rank in these metrics is mutually exclusive. That is, a course can have a high $\rho_H$ value, being highly correlated with overall mathematics GPA, without having a high $\beta_i$ value, indicating that it is highly impactful on overall mathematics GPA. Then, the fact that Complex Analysis ranks as the top course under both metrics says that it is more correlated and more impactful on overall mathematics GPA than any other mathematics course for our pure mathematics freshman admit group. 

To understand why Complex Analysis was the course with highest $\rho_H$ and $\beta_i$ values, we looked for other explanatory statistics.  We found that regular faculty more frequently teach Complex Analysis than temporary faculty do.   We believe this is one possible factor in explaining high $\rho_H$ and $\beta_i$ values for Complex Analysis. Permanent faculty have more teaching experience and are better able to assess students in their course.  As such, Complex Analysis grades more strongly correlate with the overall mathematics GPA. 
As to whether students are getting more adept training from regular faculty for the rest of their math major is beyond the scope of this work.  Moreover,  in this data set, we cannot determine  whether a regular or temporary faculty assigned a particular student's grade and so do not comment on this further.



\begin{table*}
\centering
\caption{Summary statistics for Pure Mathematics freshman admits. Complex Analysis is the top ranked course under both correlation and regression metrics. }
\begin{tabular}{m{4.5cm} ccccc}
\toprule
     \textbf{Course} &\textbf{$\rho_H$ Rank} &\textbf{$\rho_H$} &\textbf{$\beta_i$ Rank} &\textbf{$\beta_i$} \\
    \hline
\begin{filecontents*}{pm_results.csv}
Course,CorrelationRank,CorrelationMetric,RegressionRank,RegressionMetric
Complex Analysis,1,1.039,1,0.3774
Topology,2,1.038,--,0
Linear Algebra,5,0.951,--,0
Real Analysis (Part I),6,0.9507,--,0
Probability (Part I),10,0.8224,2,0.00824
Honors Abs. Algebra (Part I),11,0.811,3,0.01919
Diff. Eqs. (Lower Div.),14,0.8012,4,0.001573
\end{filecontents*}
\csvreader[ late after line=\\\hline,table foot = \bottomrule]
{pm_results.csv}{Course = \Course, CorrelationRank = \CorrelationRank, CorrelationMetric = \CorrelationMetric, RegressionRank = \RegressionRank, RegressionMetric = \RegressionMetric  }
{ \Course & \CorrelationRank & \CorrelationMetric & \RegressionRank & \RegressionMetric}
\end{tabular}
\label{table:results_final}
\end{table*}

\begin{figure}[h!]
\centering
\includegraphics[scale = 0.7]{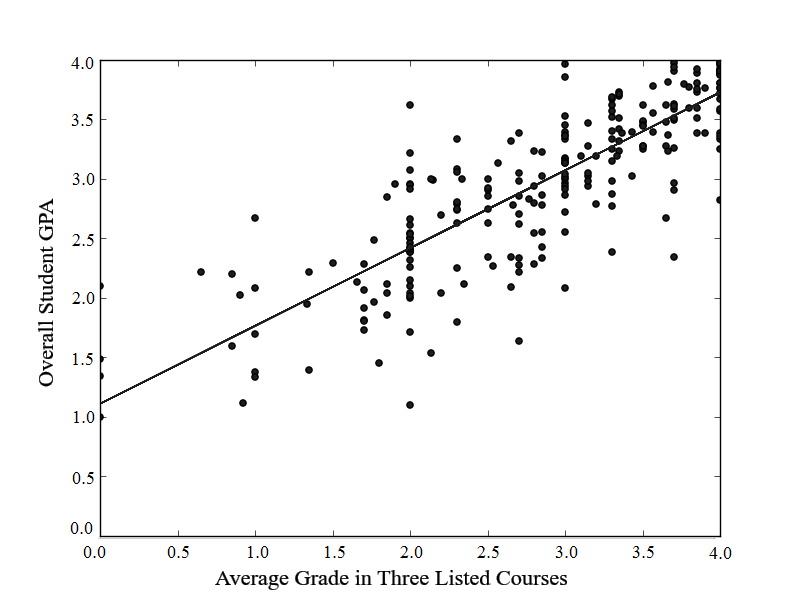}
\caption{Scatter plot of grade values in Complex Analysis, Linear Algebra, and Probability Theory (Part I) plotted against overall GPA with ridge prediction function.}
\label{fig:ridge_pure}
\end{figure}

Having selected some courses from the lasso model, we can now generate a prediction function to estimate overall mathematics GPA for pure mathematics freshman admit students. We choose to run ridge regression using Linear Algebra (Part I), Complex Analysis, and Probability Theory (Part I). We choose the latter two courses because they are assigned non-zero lasso coefficients as seen in Table \ref{table:results_final}. We choose Linear Algebra (Part I) because it is a departmental core course and has a relatively high correlation rank in Table \ref{table:results_final}. From our original student population of pure mathematics students who have taken at least one upper division course, we further sample the students who have taken all three courses. The results of our ridge regression are as follows
$${\hat{g}}^{\t{GPA}} = 0.74 + 0.29g^{132} + 0.27g^{115A} + 0.18g^{170A}.$$
We obtain an $R_{sq}$ of  87\%. Although we use grades in the above three courses to estimate overall mathematics GPA, in general students take far more courses. In this sample alone, the average number of courses taken is 16 courses with standard deviation of 4.11. Taking three standard deviations below the mean of 16, we still arrive at a number higher than 3 courses, indicating that it is unlikely to find students whose overall mathematics GPA's are determined only by the three courses used in the model. In Figure \ref{fig:ridge_pure} we can see a scatter plot of the results.


\FloatBarrier

\subsection{Survey of Other Mathematics Majors}

While we have focused primarily on the Pure Mathematics freshman admit group for our analyses, we have also used the outlined methods on two other large groups of mathematics majors. Namely, we have obtained results for the Math/Econ freshman admit group and the Applied Mathematics freshman admit group. We display the results for these groups in Table \ref{table:me_results} and Table \ref{table:am_results}. 

As we can observe from Table \ref{table:me_results}, within the Math/Econ freshman group, there are two surprising results. Firstly, the course Abstract Algebra (Part I), which is not required for Math/Econ students, has a high $\rho_H$ value. Additionally, Real Analysis (Part I), which is a core course, has a high $\beta_i$ value. The result of the course grade in Real Analysis (Part I) having high impact is surprising because of the few number of theoretical courses students of the Math/Econ major are required to take.  Additionally, the result of a course grade in Abstract Algebra having high correlation with one's overall mathematics GPA is surprising because it is both theoretical and not required within the major. Although theoretical content may prepare students for a well-rounded foundation for their future coursework, most courses within the Math/Econ curriculum tend to be applied. 

As we can observe from Table \ref{table:me_results} and Table \ref{table:am_results} for both the Math/Econ and Applied Mathematics freshman groups, Complex Analysis does not have highly ranked $\rho_H$ and $\beta_i$ values in either of these student groups. This contrasts with the results found in Section \ref{pm_study}.
Another result to highlight is that in both the Applied Math freshman group and the Math/Econ freshman group, Linear Algebra is not selected as an impactful course by lasso regression. Furthermore, in the Applied Mathematics freshman group, Real Analysis (Part I) is also not selected as an impactful course through lasso regression. 
Although both courses are considered core by the mathematics department, other courses, such as Probability Theory (Part II) and Mathematical Modeling are more correlated with and impactful on overall mathematics GPA.
\begin{table*}
\centering
\caption{Summary statistics for Math/Econ freshman admits.  }
\begin{tabular}{m{3.5cm} cccc}%
\toprule
     \textbf{Course} &\textbf{$\rho_H$ Rank} & \textbf{$\rho_H$}  &\textbf{$\beta_i$ Rank} &\textbf{$\beta_i$} \\
    \hline
\begin{filecontents*}{me_results.csv}
Course,CorrelationRank,CorrelationMetric,RegressionRank,RegressionMetric
Abs. Algebra (Part I),1,1.08271,--,0
Single Var. Calc. (Part I),2,1.02221,--,0
Multivar. Calc. (Part II),3,1.01405,3,0.68589
Num. Analysis (Part II),4,1.00849,--,0
Diff. Eq. (Low. Div.),5,0.986724,7,0.010956
Linear Algebra (Part I),6,0.94111,--,0
Real Analysis (Part I),13,0.9101,1,0.08471
Optimization,12,0.916562,2,0.081447
\end{filecontents*}
\csvreader[ late after line=\\\hline,table foot = \bottomrule]
{me_results.csv}{Course = \Course, CorrelationRank = \CorrelationRank, CorrelationMetric = \CorrelationMetric, RegressionRank = \RegressionRank, RegressionMetric = \RegressionMetric  }
{ \Course & \CorrelationRank & \CorrelationMetric & \RegressionRank & \RegressionMetric}
\end{tabular}
\label{table:me_results}
\end{table*}

\begin{table*}
\centering
\caption{Summary statistics for Applied Mathematics freshman admits. }
\begin{tabular}{m{3.5cm} cccc}
\toprule
     \textbf{Course} &\textbf{$\rho_H$ Rank} & \textbf{$\rho_H$}  &\textbf{$\beta_i$ Rank} &\textbf{$\beta_i$} \\
    \hline
\begin{filecontents*}{am_results.csv}
Course,CorrelationRank,CorrelationMetric,RegressionRank,RegressionMetric
Probability (Part II),1,1.0759,5,0.02965
Linear/Nonlinear Diff. Eq.,2,1.0714,--,0
Combinatorics,3,1.0616,--,0
Math. Modeling,4,1.055,1,0.08497
Optimization,5,1.0506,--,0
Num. Analysis (Part II),11,0.9966,2,0.06327
Complex Analysis,9,1.01803,3,0.04523
Real. Analysis (Part I),13,0.98749,--,0
Linear Algebra (Part I),15,0.95237,--,0
\end{filecontents*}
\csvreader[ late after line=\\\hline,table foot = \bottomrule]
{am_results.csv}{Course = \Course, CorrelationRank = \CorrelationRank, CorrelationMetric = \CorrelationMetric, RegressionRank = \RegressionRank, RegressionMetric = \RegressionMetric  }
{ \Course & \CorrelationRank & \CorrelationMetric & \RegressionRank & \RegressionMetric}
\end{tabular}
\label{table:am_results}
\end{table*}

\FloatBarrier

\section{Summary and Interpretation of Results}


In this paper, we developed methods to analyze the core courses selected by the mathematics department. We identified statistical attributes we expect the core courses to have. The first attribute was a correlated metric, denoted by $\rho_H$, which measured the correlation between a course grade and overall mathematics GPA. The second attribute was a metric of impact, denoted by $\beta_i$, which measured the impact a course grade has on overall mathematics GPA. We tested our methods on various mathematics majors at UCLA to validate the selection of core courses as selected by the mathematics department. 

The most surprising results came from the pure mathematics freshman group. We summarize our results for pure mathematics students in Table \ref{table:results_final}. In this table, one can see that Complex Analysis appears in both results for the mixed correlation metric and lasso regression. This suggests that the grade in Complex Analysis is both highly correlated to and impactful on one's overall mathematics GPA. Although this result is unique in that there is a single course that is highly correlated to and  highly impactful on overall mathematics GPA within this student group, it may also suggest a weakness in our methods, as other student groups did not have a single course that ranked highest by both measures.

We can interpret the results from the pure mathematics freshmen by looking at instructors throughout the 15 year period. We see that  Complex Analysis is more frequently taught by regular faculty than many other upper division courses at UCLA. In particular, about 52$\%$ of the professors who taught Complex Analysis were regular faculty. We can compare this to 46$\%$ of regular faculty who taught Real Analysis, and 39$\%$ of regular faculty who taught Linear Algebra. This may help explain the high correlation attributed to Complex Analysis. For instance, regular faculty may tend to give grades more indicative of a student's overall performance. Additionally, regular faculty may teach a course more consistently than temporary faculty so that the grade in the course is a more accurate reflection of the student's performance.  If this is the case, the math department may benefit from assigning more regular faculty to courses which are seen as foundational.


\FloatBarrier

\FloatBarrier

\subsection{Future Work}

Our core course analysis could be expanded upon with the inclusion of more data and by refining the students we consider. We purposefully selected versatile models to illustrate this core course analysis and how it could be used to better understand the mathematics curriculum. Going forward, we will utilize more directed and specific statistical models, derive the results from these models and compare them to the results presented in this paper.

We would also like to combine the lasso and ridge models to incorporate the sparsity of the $\ell_1$ norm with the metrics coming from the $\ell_2$ norm. In addition, we can expand our regression models to include explanatory variables such as average course load, whether or not the student switched majors, and how long a student took to declare a major. The goal of including such variables is to introduce more sources of variation in overall GPA into our model than just grades. We can analyze course rankings under the $\beta_i$ values from this more comprehensive regression model and compare the associated rankings to those outlined in this paper. We can then investigate any discrepancies between the two rankings and provide the findings to the department which can act accordingly. 


\appendix\section{Student Data}
\label{app1}
\label{course_appendix}

In this appendix, we briefly discuss the data that we investigated in this paper.  As mentioned earlier, each row represents a course a student took, his grade in that course, and other pertinent information.  The so-called Student Classification in Table \ref{table:data} indicates their academic standing: freshman (UFR), sophomore (USO), junior (UJR), and senior (USR).  The Admit Class is the academic standing when a student \emph{first} entered UCLA.  A first-year non-transfer student may enter as as USO or a UFR based on the number of college level credits earned while in high school. The grade values range from 0.0-4.3, where 4.3 indicates an A+ grade.  Under the umbrella of the mathematics department their are two possible subject headings: MATH and COMPT, the latter representative of the computer based courses in the program for computation.  The Course codes can be viewed at the \href{https://www.math.ucla.edu/ugrad/courses}{UCLA Department of Mathematics}.
The Major Codes and corresponding major names can be viewed at the \href{http://www.registrar.ucla.edu/dept/majorug.htm}{UCLA Registrar}.

\begin{table*}
\centering
\caption{Sample of data provided by department of mathematics. }
\begin{tabular}{m{1.6cm} cccccccc}%
\toprule
     \textbf{ID} &\textbf{Term} & \textbf{Major Code}  &\textbf{Classification} &\textbf{Admit Class} &\textbf{Subject} &\textbf{Course} &\textbf{Grade} \\
    \hline
\begin{filecontents*}{data_sample.csv}
ID,Term,Major Code,Classification,Admit Class,Subject,Course,Grade
0x3915,09W,756,USO,UFR,MATH,0003A,3.3
0xE443,09W,671,USO,UFR,MATH,0003A,3
0XC044,09W,00X,USR,UFR,COMPT,10A,3.7
\end{filecontents*}
\csvreader[ late after line=\\\hline,table foot = \bottomrule]
{data_sample.csv}{ID = \ID, Term = \Term, Major Code = \MajorCode, Classification = \Classification, Admit Class = \AdmitClass, Subject = \Subject, Course = \Course, Grade = \Grade}
{ \ID & \Term & \MajorCode & \Classification & \AdmitClass & \Subject & \Course & \Grade}
\end{tabular}
\label{table:data}
\end{table*}

\section{Acknowledgements}
We would like to thank Dimitri Shlyaktenko and Charles Z. Marshak for all the guidance they have provided us during this project and the writing of this paper. 
We thank Dimitri Shlyakhtenko and Andrea Bertozzi for their help with administrative issues related to use of this dataset. We also want to thank Mihai Cucuringu for initially helping us parse through the data and M. Puck Rombach for additional oversight and helpful discussions. Finally, we would like to thank our REU group members Dillon Montag, Milica Had\^zi-Tanovi\'c, and Juan Carlos Apitz.  
Ritvik Kharkar was supported by NSF grant DMS-1045536 and Jessica Tran was supported by the UC LEADS program.

\bibliographystyle{plain}
\bibliography{bibliography}
\FloatBarrier




\end{document}